\documentclass[letterpaper, 10 pt, conference]{ieeeconf}  

\IEEEoverridecommandlockouts 

\usepackage{afterpage}
\usepackage{epsfig} 
\usepackage{times} 
\usepackage{amsmath} 
\usepackage{amssymb}  
\usepackage{amsfonts}
\usepackage{mathptmx} 
\usepackage{tabularx}
\usepackage{epstopdf}
\usepackage{float}
\usepackage{epic,color}
\usepackage{mathrsfs}
\usepackage{dsfont,MnSymbol}
\usepackage{algorithm}
\usepackage{multirow} 
\usepackage[T1]{fontenc}

\newcounter{rmnum}

\newcounter{anum}

\def\IEEEQEDclosed{\mbox{\rule[0pt]{1.3ex}{1.3ex}}}

\def\qed{\ifmmode\IEEEQEDclosed\else{\unskip\nobreak\hfil
		\penalty50\hskip1em\null\nobreak\hfil\IEEEQEDclosed
		\parfillskip=0pt\finalhyphendemerits=0\endgraf}\fi}

\def\qed{\hspace*{\fill}~\IEEEQED\par\endtrivlist\unskip}

\def\Re{\mathbb{R}}

\def\notes#1{\marginpar{\tiny #1}\typeout{Notes!
Notes!
Notes!
}}
\renewcommand{\notes}[1]{\typeout{notes!}}
\def\FRAC#1#2#3{\genfrac{}{}{}{#1}{#2}{#3}}
\def\half{{\mathchoice{\FRAC{1}{1}{2}}%
{\FRAC{2}{1}{2}}%
{\FRAC{3}{1}{2}}%
{\FRAC{4}{1}{2}}}}

\def\Re{\field{R}}

\def\clP{{\cal P}}
\def\clZ{{\cal Z}}

\def\E{{\sf E}}

\def\IEEEQEDclosed{\mbox{\rule[0pt]{1.3ex}{1.3ex}}}
\def\qed{\nobreak\hfill\IEEEQEDclosed}

\def\clZ{{\cal Z}}

\newtheorem{theorem}{Theorem}

\newtheorem{definition}{Definition}

\newtheorem{remark}{Remark}
\newtheorem{proposition}{Proposition}

\def\beq{\begin{eqnarray}} 
\def\bc{\begin{center}} 
\def\be{\begin{enumerate}}
\def\bi{\begin{itemize}} 
\def\bs{\begin{small}}
\def\bS{\begin{slide}}
\def\ec{\end{center}} 
\def\ee{\end{enumerate}}
\def\ei{\end{itemize}}
\def\es{\end{small}}
\def\eS{\end{slide}}
\def\eeq{\end{eqnarray}}

\newcommand{\newP}[1]{\medskip\noindent{\bf #1:}}

\newcommand{\ud}{\,\mathrm{d}}

\def\Re{\mathbb{R}}
\def\E{{\sf E}}

\def\clP{{\cal P}}
\def\clZ{{\cal Z}}

\renewcommand{\Re}{\mathbb{R}}

\def\FRAC#1#2#3{\genfrac{}{}{}{#1}{#2}{#3}}

\def\sJ{{\sf J}}
\def\sP{{\sf P}}
\def\tsP{{\tilde{\sf P}}}

\def\clF{{\cal F}}

\def\clP{{\cal P}}
\def\clU{{\cal U}}
\def\clV{{\cal V}}

\def\clZ{{\cal Z}}
\def\E{{\sf E}}

\def\tE{\tilde{\sf E}}

\def\bS{\mathbb{S}}

\newlength{\noteWidth}
\long\def\notes#1{\ifinner
	{\tiny #1}
	\else
	\marginpar{\parbox[t]{\noteWidth}{\raggedright\tiny #1}}
	\fi}

\title{\LARGE \bf A Dynamic Programming Formulation for the Nonlinear Filter}

\author{Jin Won Kim and Prashant G. Mehta
	\thanks{Financial support from the 
		NSF grant 1761622 and the 
		ARO grant W911NF1810334 is gratefully acknowledged. 
	}
	\thanks{J.~W. Kim and P.~G.~Mehta are with the Coordinated
		Science Laboratory and the Department of Mechanical Science and
		Engineering at the University of Illinois at Urbana-Champaign
		(UIUC); Corresponding email: mehtapg@illinois.edu.}
}

\begin{document}

\maketitle
\thispagestyle{empty}
\pagestyle{empty}

\begin{abstract}
	
This paper build on our recent work where we presented a dual stochastic optimal control formulation of the nonlinear filtering problem~\cite{kim2019duality}.  The constraint for the dual problem is a backward stochastic differential equations (BSDE).  The solution is obtained via an application of the maximum principle (MP).  In the present paper, a dynamic programming (DP) principle is presented for a special class of BSDE-constrained stochastic optimal control problems.  The principle is applied to derive the solution of the nonlinear filtering problem. 


\end{abstract}

\section{Introduction}
\label{sec:intro}

A backward stochastic differential equation (BSDE) is an It\^o
stochastic differential equation over a finite time-horizon $[0,T]$
where the terminal condition at the terminal time $T$ is
specified. BSDE was first introduced by Bismut as an adjoint equation
for the linear-quadratic stochastic optimal control
problem~\cite{bismut1978introductory}. Later, Pardoux and
Peng~\cite{pardoux1990adapted} introduced nonlinear BSDEs and proved
the existence and uniqueness results for the general Lipschitz
cases. BSDEs have several applications in  stochastic optimal
control~\cite{kohlmann2000relationship}, mathematical
finance~\cite{el1997backward}, and analysis of certain types of
partial differential equations~\cite{pardoux1998backward}.
Since BSDE is a dynamical system, it is natural to investigate
control problems for BSDEs. After Peng~\cite{peng1993backward} introduced the optimal control
problem for forward-backward SDEs, the theoretical framework and
applications for the stochastic maximum principle (MP) is widely
studied~\cite{dokuchaev1999stochastic,lim2001linear,ji2006maximum,huang2009maximum,wang2018introduction}.

In a recent paper~\cite{kim2019duality}, we introduced a dual model to
transform the classical nonlinear filtering problem into a
BSDE-constrained stochastic optimal control problem.  The model has
since been used for the purposes of defining observability of the nonlinear filtering problem~\cite{kim2019observability}, and for analysis of the filter stability in the ergodic~\cite{kim2021ergodic} and non-ergodic~\cite{kim2021detectable} settings of the problem.
In the original paper~\cite{kim2019duality}, the solution of the optimal control problem is obtained by using the MP.

The objective of the present paper is to introduce a dynamic programming (DP)
approach to solve a class of BSDE-constrained stochastic optimal
control problems.  After the DP equation is presented in
Section~\ref{sec:main}, it is applied to the dual optimal control
problem of nonlinear filtering.

The outline of the remainder of this paper is as follows: The
nonlinear filtering problem and its dual -- the BSDE-constrained
optimal control problem -- appear in Sec.~\ref{sec:duality}. The DP
equation and its application to the nonlinear filtering problem are
presented in~Sec.~\ref{sec:main}. 
The Appendix contains the proofs.

\section{Duality for nonlinear filtering}\label{sec:duality}



\newP{Notation} The state-space $\bS:=\{1,2,\hdots,d\}$ is finite.  
The set of probability vectors on $\bS$ is denoted by $\clP(\bS)$:
$\mu\in \clP(\bS)$ if $\mu(x)\geq 0$ and $\sum_{x\in\bS} \mu(x) =
1$. The space of deterministic functions on $\bS$ is identified with
$\Re^d$: Any function $f:\mathbb{S}\to \Re$ is determined by its
value $f(x)$ at $x\in \bS$.  For a measure $\mu \in {\cal
	P}(\mathbb{S})$ and a function $f\in \Re^d$,   $\mu(f):=\sum_x
\mu(x) f(x)$. For two vectors $f,h\in\Re^d$, $f\, h$ denotes the
element-wise (Hadamard) product: $(f\,h)(x) := f(x) h(x)$ and
similarly $f^2  = f\, f$.  

\subsection{Filtering problem} 

Consider a pair of continuous-time stochastic processes $(X,Z)$
defined on a probability space $(\Omega,\clF,{\sf P})$:

(1) The state $X= \{X_t\in\bS:t\geq 0\}$ is a Markov process with
initial condition $X_0\sim\mu\in{\clP}(\bS)$ (prior) and the rate matrix 
$A\in\Re^{d\times d}$.   
For a function $f\in \Re^d$, the carr\'e du
champ operator $\Gamma:\Re^d\to\Re^d$ is defined according to 
\[\Gamma (f)(x) := \sum_{j \in \mathbb{S}} A(x,j)
(f(x) - f(j))^2 \quad \text{for} \; x\in\bS 
\] 

(2) The observation process $Z = \{Z_t\in\Re:t\geq 0\}$ is defined
according to the following model:
\[
Z_t := \int_0^t h(X_\tau) \ud \tau + W_t
\]
where $h : \bS \to \Re$ is the observation function and 
$W=\{W_t\in\Re:t\geq 0\}$ is a Wiener process (w.p.) that is assumed
to be independent of $X$. The scalar-valued observation model is considered for notational ease. 
The filtration generated by $(X,W)$ is denoted by $\clF := \{\clF_t:t \ge 0\}$, and the filtration generated by $Z$ is denoted as $\clZ := \{\clZ_t : t\ge 0\}$ where $\clZ_t = \sigma(\{Z_s : s\le t\})$.

The filtering problem is to compute the conditional distribution (posterior), denoted by $\pi_t\in\clP(\bS)$, of the state $X_t$ given $\clZ_t$. For $f\in \Re^d$, $\pi_t(f) := \E\big(f(X_t)\mid \clZ_t\big)$ is the object of interest.


A standard approach~\cite[Ch.~5]{xiong2008introduction} on optimal filtering is to consider a new probability measure $\tsP$ on $\Omega$ such that the Radon-Nikodym derivative with respect to $\sP$ is given by
\[
\frac{\ud \tsP}{\ud \sP} \Big|_{\clF_t} = \exp\Big(-\int_0^t h(X_\tau )\ud W_\tau  - \half \int_0^t |h(X_\tau )|^2 \ud \tau\Big)=:D_t^{-1}
\]
Then $Z$ is a $\tsP$-Brownian motion independent of $X$. Under this new measure, the unnormalized filter is defined by $\sigma_t(f) := \tE(D_tf(X_t)|\clZ_t)$ where $\tE$ denotes the expectation with respect to $\tsP$.

\newP{Function spaces}
Under $\tsP$, $Z$ is a Brownian motion, and therefore the following Hilbert spaces are defined:
The space $L^2_\clZ([0,T];S)$ is a $\clZ$-adapted random processes taking values in $S$. Similarly, $L^2_{\clZ_T}(S)$ is $S$-valued random element which is $\clZ_T$-measurable. In both cases, the $L^2$ norm is with respect to $\tsP$. 

%

\subsection{Dual optimal control problem}
In our recent work~\cite{kim2019duality}, a dual optimal control
problem is introduced.  It is based upon the following linear BSDE:
\begin{equation}\label{eq:dual-bsde}
\begin{aligned}
- \ud {Y}_t(x) & = \big((A {Y}_t)(x) + h(x)\big(U_t + 
V_t(x)\big)\big)\ud t - {V}_t(x)\ud Z_t, \\
Y_T (x) &= F(x) \;\; \forall \; x \in \bS 
\end{aligned}
\end{equation}
The boundary condition prescribed at the terminal time $T$ is allowed to be random, with $F\in L^2_{\clZ_T}(\Re^d)$. The control $U = \{U_t:0\le t \le T\}$ is chosen in the set of admissible control $L^2_\clZ([0,T];\Re)$.
The solution pair $(Y,V) =\{(Y_t,V_t):0\le t\le T\}$ of the BSDE is adapted to the filtration $\clZ$, and it is uniquely determined in $L_\clZ^2([0,T];\Re^d) \times L_\clZ^2([0,T];\Re^d)$~\cite{pardoux1990adapted}.

Define the cost functional
\begin{equation}
\sJ(U) := \tE\Big(\half |Y_0(X_0)-\mu(Y_0)|^2 + \int_0^T l(Y_t,V_t,U_t,t)\ud t
\Big)\label{eq:dual-oc-problem}
\end{equation}
where the Lagrangian $l$ is given by
\[
l(y,v,u,t;\omega) := \half \sigma_t\big(\Gamma(y)(\cdot)\big) + \half \sigma_t\big(|u + v(\cdot)|^2\big)
\]
The {\em dual optimal control problem} is to choose a 
control $U\in L^2_\clZ([0,T];\Re)$ such that $\sJ(U)$ is minimized subject to the BSDE constraint~\eqref{eq:dual-bsde}.


The following proposition is a version of the main result in~\cite{kim2019duality}. The justification appears in Appendix~\ref{apdx:pf-duality}.

\medskip

\begin{proposition}\label{prop:duality} 
Consider the dual optimal control problem. Define
\[
S_t = \mu(Y_0)- \int_0^t U_\tau  \ud Z_\tau ,\quad t\in[0,T]
\]
Then for all $t\in[0,T]$,
\begin{equation*}\label{eq:cost-to-go-clue}
\sJ(U) = \tE\Big(\half D_t |Y_t(X_t)-S_t|^2 + \int_t^T l(Y_\tau ,V_\tau ,U_\tau ,\tau)\ud \tau\Big)
\end{equation*}
In particular,
\begin{equation}\label{eq:duality-principle}
\sJ(U) = \half \tE\big(D_T|F(X_T)-S_T|^2\big) = \half \E\big(|F(X_T)-S_T|^2\big)
\end{equation}
\end{proposition}

\medskip

The equation~\eqref{eq:duality-principle} is called \emph{duality principle}, and it connects the minimum variance estimate problem and the BSDE-constrained stochastic optimal control problem.

\section{Main results}\label{sec:main}

\subsection{Optimal control problem on BSDE}

Consider the BSDE-constrained stochastic optimal control problem in more general form.
\begin{align}
&\mathop{\mathrm{Minimize}}_{U\in{\cal U}}: \sJ_T(U;\xi):= \E\Big(h(Y_0) + \int_0^T l(Y_t, V_t, U_t,t) \ud t \Big) \label{eq:cost}\\
&\text{Subject to}: \ud Y_t = f(Y_t,V_t,U_t,t) \ud t + V_t \ud Z_t,\quad Y_T = \xi \label{eq:dynamics}
\end{align}
where $\xi \in L_{\clZ_T}^2(\Re^d)$. The set of admissible control is $ \clU = L_\clZ^2([0,T];\Re)$. The drift term $f:\Re^d\times \Re^{d} \times \Re \times [0,T]\times \Omega \to \Re^d$ is assumed to be uniformly Lipschitz with each argument almost surely, for almost every $t$. For every value of $y,v$ and $u$, $f(y,v,u,\cdot)$ is $\clZ$-adapted. Under this condition,~\eqref{eq:dynamics} admits the unique solution pair $(Y,V) \in L_\clZ^2([0,T];\Re^d) \times L_\clZ^2([0,T];\Re^{d})$~\cite{pardoux1990adapted}.
For the cost functional, we assume $h:\Re^d\to \Re$ and $l:\Re^d\times \Re^{d} \times \Re \times [0,T]\times \Omega \to \Re^d$ are bounded almost surely. Similar to $f$, $l(y,v,u,\cdot)$ is $\clZ$-adapted.


\subsection{Value function of stochastic optimal control problems}

To formulate the dynamic programming principle, consider a partial problem up to time $t \le T$ from  $\zeta\in L^2_{\clZ_t}(\Re^d)$ defined by:
\[
\sJ_t(U;\zeta) = \E\Big(h(Y_0) + \int_0^t l(Y_\tau^{\zeta,t},V_\tau^{\zeta,t},U_\tau,\tau) \ud \tau \Big)
\]
where $\big\{\big(Y_\tau^{\zeta,t}, V_\tau^{\zeta,t}\big) : \tau \in [0,t]\big\}$ is the solution to the BSDE:
\[
\ud Y_\tau^{\zeta,t} = f\big(Y_\tau^{\zeta,t},V_\tau^{\zeta,t},U_\tau,\tau\big)\ud \tau + V_\tau^{\zeta,t} \ud W_\tau,\quad Y_t^{\zeta,t} = \zeta
\]
\begin{definition}\label{def:value-function}
Consider the optimal control problem~\eqref{eq:cost}-\eqref{eq:dynamics}. The \emph{value function} is a sequence of functions $\clV:=\{\clV_t:0\le t\le T\}$ where $\clV_t:L^2_{\clZ_t} \to \Re$ is defined by:
\begin{equation}\label{eq:value-function}
\clV_t(\zeta) = \inf_{U\in{\cal U}} \sJ_t(U;\zeta)
\end{equation}
\end{definition}

\medskip

Analogously with the forward-in-time stochastic DP principle, the following theorem whose proof appears in Appendix~\ref{apdx:pf-DP-principle} is proposed:

\medskip

\begin{theorem}\label{thm:DP-principle}
Let $\clV$ be the value function of the optimal control problem~\eqref{eq:cost}-\eqref{eq:dynamics}. Then it satisfies the following:
\begin{enumerate}
	\item[(i)] $\clV_0(\cdot) = h$
	\item[(ii)] For any $0\le s < t \le T$ and any $\zeta\in L^2_{\clZ_t}(\Re^d)$,
	\begin{equation}\label{eq:DP-principle}
	\clV_t(\zeta) = \inf_{U\in{\cal U}} \E\Big(\clV_s(Y_s^{\zeta,t}) + \int_s^t l(Y_\tau^{\zeta,t},V_\tau^{\zeta,t},U_\tau,\tau)\ud \tau\Big)
	\end{equation}
\end{enumerate}
\end{theorem}

%

\medskip

For stochastic optimal control problems, an appealing formulation for the value function is to construct a martingale associated with it. The martingale version of DP principle is as follows. The proof appears in Appendix~\ref{apdx:pf-martingale-DP}.

\medskip

\begin{proposition}\label{thm:martingale-DP}
Let $\clV$ be the value function of the optimal control problem~\eqref{eq:cost}-\eqref{eq:dynamics}. Then
\begin{enumerate}
	\item[(i)] $\clV_0(\cdot) = h$
	
	\item[(ii)] Define $M^U = \{M_t^U:0\leq t\leq T\}$ for any admissible control $U\in{\cal U}$ by:
	\begin{equation}\label{eq:martingale-def}
	M_t^U = \clV_t(Y_t) - \int_0^t l(Y_\tau,V_\tau,U_\tau,\tau) \ud \tau
	\end{equation}
	where $(Y,V)$ is the solution to~\eqref{eq:dynamics}. $M^U$ is a super-martingale for any admissible control $U$; and it is a margingale if and only if $U$ is the optimal solution.
\end{enumerate}
\end{proposition}

\medskip


\subsection{Optimal control obtained via the martingale DP principle}

The following theorem whose proof appears in Appendix~\ref{apdx:pf-converse-DP} characterizes the optimal control using dynamic programming.

\medskip

\begin{theorem}\label{thm:converse-DP}
	Suppose there exists $\clV$ and $U^* \in {\cal U}$ such that:
	\begin{enumerate}
		\item[(i)] $\clV_0(\cdot) = h(\cdot)$.
		
		\item[(ii)] The process $M^U$ defined by~\eqref{eq:martingale-def} is a super-martingale for each admissible control $U\in \clU$, and a martingale for $U=U^*$.
	\end{enumerate}
	Then $U = U^*$ is an optimal control with cost $\E\big(\clV_T(\xi)\big)$.
\end{theorem}

\medskip

\begin{remark}
	From the definition of super-martingale, the second condition is equivalent to write for any $0\le s\le t\le T$,
	\begin{equation}\label{eq:Martingale-DP-conclusion}
	\clV_s(Y_s) \geq \E\Big(\clV_t(Y_t) - \int_s^t l(Y_\tau,V_\tau,U_\tau,\tau) \ud \tau\; \big\mid\, \clZ_s \Big)
	\end{equation}
	For forward-in-time Markovian stochastic control problems, the counterpart of~\eqref{eq:Martingale-DP-conclusion} is exactly the dynamic programming principle and $\clV$ is the \emph{value function} (cf.~\cite[Remark 6.1.5]{van2007stochastic}).
	However for BSDE problems, conditioning on $\clZ_s$ is not the
	same as fixing on $Y_s$. Therefore, the conditions in
	Theorem~\ref{thm:converse-DP} do not yield an interpretation of $\clV_t$ as the value function at time $t$ in this case but only concludes that $U^*$ is optimal.
\end{remark}

\subsection{Application to nonlinear filtering}

The Proposition~\ref{prop:duality} suggests that the optimal solution
yields $S_t = \pi_t(Y_t)$. Hence, consider $\clV_t(\zeta)$ to be
\begin{equation}\label{eq:value-function-filter}
\clV_t(\zeta) := \half \tE\big(D_t|\zeta(X_t) - \pi_t(\zeta)|^2 \big)
\end{equation}
The following proposition whose proof appears in Appendix~\ref{apdx:pf-duality-DP} allows to obtain the optimal solution in a DP approach.

\medskip

\begin{proposition}\label{prop:duality-DP}
Consider the dual optimal control problem~\eqref{eq:dual-bsde}-\eqref{eq:dual-oc-problem}. Let $M^U$ be defined by~\eqref{eq:martingale-def} where $\clV_t$ is defined as~\eqref{eq:value-function-filter}. Then $M^U$ is a $\tsP$-super-martingale, and $M^U$ is a $\tsP$-martingale if and only if for all $t$,
\begin{equation}\label{eq:optimal-solution}
U_t = -\big(\pi_t(h Y_t) - \pi_t(h)\pi_t(Y_t)\big) - \pi_t(V_t)
\end{equation}
\end{proposition}

\medskip

\begin{remark}
By the Theorem~\ref{thm:converse-DP}, we
conclude~\eqref{eq:optimal-solution} is the optimal solution to the
dual optimal control problem. In~\cite{kim2019duality}, stochastic
maximum principle is used to derive the optimal control. A similar
super-martingale is also considered using the innovation process.  The
change of measure is introduced to make $\clZ$ be a filtration
generated by a Brownian motion $Z$ in this paper.
\end{remark}

\section{Conclusion}\label{sec:conclusion}

The DP approach provides a sufficient condition to obtain the optimal control directly, while MP provides necessary condition.
However, a verification theorem to obtain the value function for control problems on BSDEs is still an open question. Although the function $\clV$ plays a role like the value function in its forward-in-time counterpart, it is challenging to prove that $\clV$ is in fact the value function due to the information structure. 
This is subject to future research.


\bibliographystyle{IEEEtran}
\bibliography{_master_bib_jin,jin_papers,bsde-ctrl}

\appendix
\section{Appendix}

\subsection{Proof of the Proposition~\ref{prop:duality}}\label{apdx:pf-duality}

The claim is essentially~\cite[Prop.~1 and 2]{kim2019duality}. Therefore the only justification is that the optimal control formulation is identical.
The original formulation in~\cite{kim2019duality} is:
\begin{align*}
\sJ(U) &= \E\Big(\half|Y_0(X_0)-\mu(Y_0)|^2 \\
&\qquad+ \int_0^T \half \Gamma (Y_t)(X_t) +\half|U_t+V_t(X_t)|^2 \ud t\Big)
\end{align*}
Apply the change of measure from $\sP$ to $\tsP$:
\begin{align*}
\sJ(U) &= \tE\Big(\half D_0|Y_0(X_0)-\mu(Y_0)|^2 \\
&\qquad+ \int_0^T \half D_t\Gamma (Y_t)(X_t) +\half D_t|U_t+V_t(X_t)|^2 \ud t\Big)
\end{align*}
By the tower property of conditional expectation,
\begin{align*}
\sJ(U) &= \tE\Big(\half|Y_0(X_0)-\mu(Y_0)|^2 \\
&\qquad + \int_0^T \half \tE\big(D_t \Gamma(Y_t)(X_t)\mid \clZ_t\big) \ud t\\
&\qquad +\int_0^T  \half \tE(D_t|U+V_t(X_t)|^2\mid \clZ_t) \ud t\Big)\\
&= \tE\Big(\half|Y_0(X_0)-\mu(Y_0)|^2  + \int_0^T l(Y_t,V_t,U_t,t)\ud t \Big)
\end{align*}
Therefore,~\eqref{eq:dual-oc-problem} is identical to the cost functional considered in~\cite{kim2019duality}, and the claim follows.

\medskip

\subsection{Proof of the Theorem~\ref{thm:DP-principle}}\label{apdx:pf-DP-principle}

We start from the definition of the value function:
\begin{align*}
\clV_t(\zeta) &= \inf_{U\in{\cal U}}  \E\Big(h(Y_0^{\zeta,t}) + \int_0^t l(Y_\tau^{\zeta,t},V_\tau^{\zeta,t},U_\tau,\tau) \ud \tau \Big)\\
&= \inf_{U\in{\cal U}}  \E\Big(h(Y_0^{\zeta,t}) + \int_0^s l(Y_\tau^{\zeta,t},V_\tau^{\zeta,t},U_\tau,\tau) \ud \tau \\
&\qquad\qquad\qquad+ \int_\tau ^t l(Y_\tau^{\zeta,t},V_\tau^{\zeta,t},U_\tau,\tau) \ud \tau \Big)
\end{align*}
The claim is that the first two terms are precisely $\clV_s(Y_s^{\zeta,t})$. Recall the integral formulae:
\[
Y_s^{\zeta,t} = \zeta - \int_s^t f(Y_\tau,V_\tau,U_\tau,\tau) \ud \tau - \int_s^t V_\tau \ud Z_\tau
\] 
Note that it depends only on $U_\tau: \tau\in[s,t]$. Meanwhile,
\[
Y_u^{\zeta,t} = Y_s^{\zeta,t} - \int_u^t f(Y_\tau,V_\tau,U_\tau,\tau) \ud \tau - \int_u^s V_\tau \ud Z_\tau,\quad u\le s
\]
depends only on $U_\tau:\tau\in[0,s]$ given $Y_s^{\zeta,t}$, and therefore \[
\big(Y_u^{\zeta,t},V_u^{\zeta,t}\big) = \Big(Y_u^{Y_s^{\zeta,t},s},V_u^{Y_s^{\zeta,t},s}\Big)
\]

\medskip

\subsection{Proof of the Proposition~\ref{thm:martingale-DP}}\label{apdx:pf-martingale-DP}
We start from~\eqref{eq:DP-principle}:
\[
\clV_t(\zeta) \leq \E\Big(\clV_s(Y_s^{\zeta,t}) + \int_s^t l(Y_\tau^{\zeta,t},V_\tau^{\zeta,t},U_\tau,\tau)\ud \tau \Big)
\]
Note that both sides are map a random variable $\zeta$ to a scalar. For $\zeta = Y_t$, 
\[
(Y_s,V_s) = \big(Y_s^{\zeta,t},V_s^{\zeta,t}\big)
\]
and therefore
\[
\E\big(\clV_t(Y_t)\mid \clZ_s \big) \leq \E\Big(\clV_s(Y_s) + \int_s^t l(Y_\tau,V_\tau,U_\tau,\tau)\ud \tau \big\mid \clZ_s \Big)
\]
Upon subtracting $\E\big(\int_0^t l(Y_\tau,V_\tau,U_\tau,\tau)\ud \tau \mid \clZ_s \big)$ on both sides, we have
\begin{align*}
\E\Big(\clV_t(Y_t) - \int_0^t &(Y_\tau,V_\tau,U_\tau,\tau)\ud \tau \,\big\mid\, \clZ_s \Big) \\
&\leq \E\Big(\clV_s(Y_s) - \int_0^s l(Y_\tau,V_\tau,U_\tau,\tau)\ud \tau \, \big\mid\, \clZ_s \Big)
\end{align*}
Since the right-hand side is $\clZ_s$-measurable, we may drop conditional expectation, and hence
$$
\E\big(M_t^U\mid \clZ_s\big)\leq M_s^U
$$
Therefore, $M^U$ is a super-martingale. The inequality becomes equality upon choosing the optimal control.

\medskip

\subsection{Proof of the Theorem~\ref{thm:converse-DP}}\label{apdx:pf-converse-DP}

By assumption that $M^U$ is a super-martingale,
$$\E \big(M_T^U\big) \leq M_0^U = \clV_0(Y_0) = h(Y_0)$$
Take expectation on the right-hand side and expand the left-hand side as
$$
\E\Big(\clV_T(\xi) - \int_0^T l(Y_t,V_t,U_t,t)\ud t\Big) \leq \E\big(h(Y_0)\big)
$$
Therefore we have
$$
\E\big(\clV_T(\xi)\big) \leq \sJ(U),\quad \forall U
$$
where equality holds for $U=U^*$.

\medskip

\subsection{Proof of the Proposition~\ref{prop:duality-DP}}\label{apdx:pf-duality-DP}
\def\bsig{{\bar{\sigma}}}

By the tower property of the conditional expectation,
\[
\clV_t(Y_t) = \half\tE\Big( \tE\big(D_t|Y_t(X_t)-\pi_t(Y_t)|^2\mid \clZ_t\big)\Big)
\]
The term inside the expectation equals to
\begin{align*}
\E\big(|Y_t(X_t)&-\pi_t(Y_t)|^2\mid\clZ_t\big)\\
=&  \tE\big(D_t\big(Y_t^2(X_t) - 2Y_t(X_t)\pi_t(Y_t) + (\pi_t(Y_t))^2\big)\mid\clZ_t\big)\\
=& \sigma_t(Y_t^2) - \sigma_t(Y_t)\pi_t(Y_t)
\end{align*}
where we used $\sigma_t(Y_t) = \sigma_t(1)\pi_t(Y_t)$.

The first term requires $\ud \big(Y_t^2\big)$, which is
\[
\ud Y_t^2 = -2Y_t(AY_t+h(U_t+V_t))\ud t + V_t^2\ud t + 2Y_tV_t \ud Z_t
\]
Use Zakai equation~\cite[Theorem 5.5]{xiong2008introduction} to take differential form of each term:
\begin{align*}
\ud \sigma_t(Y_t^2) &= \sigma_t(A Y_t^2) \ud t + \sigma_t(hY_t^2)\ud Z_t\\
&\quad +\sigma_t\big(-2Y_t\big(A Y_t + h (U_t+V_t)\big) + V_t^2\big)\ud t \\
&\quad+ 2\sigma_t\big(Y_tV_t\big) \ud Z_t +2\sigma_t(h Y_tV_t)\ud t\\
&=\sigma_t\big(\Gamma Y_t\big) \ud t + \sigma_t(V_t^2)\ud t - 2\sigma_t(hY_t)U_t \ud t\\
&\quad +\big(\sigma_t(hY_t^2)+ 2\sigma_t(Y_tV_t)\big)\ud Z_t
\end{align*}
Again use Zakai equation to compute:
\begin{align*}
\ud \sigma_t(Y_t) &= \sigma_t(AY_t)\ud t + \sigma_t(hY_t)\ud Z_t -\sigma_t\big(AY_t+h(U_t+V_t)\big)\ud t \\
&\quad + \sigma_t(V_t)\ud Z_t + \sigma_t(hV_t)\ud t\\
&=-\sigma_t(h)U_t \ud t + \big(\sigma_t(h Y_t)+\sigma_t(V_t)\big)\ud Z_t
\end{align*}
Upon using the nonlinear filter in a similar way, one can obtain
\begin{align*}
\ud \pi_t(Y_t)
&=-\pi_t(h)(U_t - U_t^*)\ud t + U_t^*\ud Z_t
\end{align*}
where $U_t^* = -\big(\pi_t(h Y_t) - \pi_t(h)\pi_t(Y_t)\big) - \pi_t(V_t)$ is as given in~\eqref{eq:optimal-solution}.
Therefore by It\^o product rule,
\begin{align*}
\ud \sigma_t(Y_t)\pi_t(Y_t) &= -\sigma_t(h)U_t\pi_t(Y_t) \ud t - \sigma_t(Y_t) \pi_t(h)(U_t-U_t^*) \ud t \\
&\quad+ \big(\sigma_t(h Y_t) + \sigma_t(V_t)\big)U_t^* \ud t + (\cdots)\ud Z_t
\end{align*}
where the martingale term is omitted.
$\ud M_t^U$ can now be simplified by collecting terms, and
\[
\ud M_t^U = -\half \sigma_t(1)(U_t-U_t^*)^2 \ud t + (\cdots)\ud Z_t
\]
Since $-\half \sigma_t(1)(U_t - U_t^*)^2 \le 0$ and $Z$ is a $\tsP^\mu$-martingale, $M^U$ is a $\tsP$-super-martingale, and it is a martingale if and only if $U_t = U_t^*$ for all $t$.

\end{document}